\documentclass{svmult}
\usepackage{etex}
\usepackage{mathptmx}
\usepackage{helvet}
\usepackage{courier}
\usepackage{makeidx}
\usepackage{graphicx}
\usepackage{multicol}
\usepackage{footmisc}
\usepackage{amsmath}
\usepackage{amssymb}
\usepackage{mathtools}
\usepackage{empheq}
\usepackage{booktabs}
\usepackage{mathtools}
\usepackage{multirow}

\usepackage{pgfplots}
\pgfplotsset{compat=newest}
\newlength\fwidth
\newlength\fheight

\newcommand{\R}{\ensuremath{\mathbb R}}
\newcommand{\U}{\ensuremath{\mathbb U}}
\newcommand{\yl}{{\ensuremath{y^\ell}}}

\DeclareMathOperator*{\argmin}{arg\,min}

\newcommand{\Wpod}{\ensuremath{W_{\text{POD}}}}
\newcommand{\Vpod}{\ensuremath{V_{\text{POD}}}}
\newcommand{\Wpodadj}{\ensuremath{W_{\text{POD,adj}}}}

\newcommand{\Wbt}{\ensuremath{W_{\text{BT}}}}
\newcommand{\Vbt}{\ensuremath{V_{\text{BT}}}}
\newcommand{\Wricc}{\ensuremath{W_{\text{Ricc}}}}

\newcommand{\cV}{\ensuremath{\mathcal V}}

\begin{document}
\title*{Model order reduction approaches for\\ infinite horizon optimal control \\problems via the HJB equation}
\titlerunning{MOR approaches for inifinite horizon OCP via the HJB equation}
\author{A. Alla \and A. Schmidt \and B. Haasdonk}
\institute{Alessandro Alla \at Florida State University, Department of Scientific Computing,\\Tallahassee, USA, \email{aalla@fsu.edu}
\and Andreas Schmidt \at University of Stuttgart, Institute for Applied Analysis and Numerical Simulation,\\ Stuttgart, Germany, \email{schmidta@mathematik.uni-stuttgart.de}
\and Bernard Haasdonk \at University of Stuttgart, Institute for Applied Analysis and Numerical Simulation,\\ Stuttgart, Germany, \email{haasdonk@mathematik.uni-stuttgart.de}}

\maketitle

\abstract{We investigate feedback control for infinite horizon
  optimal control problems for partial differential equations. The
  method is based on the coupling between Hamilton-Jacobi-Bellman
  (HJB) equations and model reduction techniques. It is well-known that
  HJB equations suffer the so called {\em curse of dimensionality} and, therefore, a
  reduction of the dimension of the system is mandatory. In this
  report we focus on the infinite horizon optimal control problem with quadratic cost functionals.
  We compare several model
  reduction methods such as Proper Orthogonal Decomposition, Balanced
  Truncation and a new algebraic Riccati equation based approach. Finally, we present
  numerical examples and discuss several features of the different
  methods analyzing advantages and disadvantages of the reduction
  methods.}

\section{Introduction}
\label{sec:introduction}
The approximation of optimal control problems for partial differential equations (PDEs) is a very challenging topic. Although it has been successfully studied for open-loop problems (we address the interested reader to the books \cite{HPUU09,Tro10} for more details), the closed-loop control problem presents several open questions for infinite dimensional equations. \\
One common way to obtain a feedback control is by means of the dynamic programming principle (DPP). The DPP characterizes the value function and its continuous version leads to a HJB equation. The theory of the viscosity solution allows us to characterize the value function as the unique solution of the HJB equations. We note that these results are quite general and valid for any problem dimension. We refer to the book \cite{BC97} for more details about the topic for ordinary differential equations. For the sake of completeness, we also mention Model Predictive Control as alternative to obtain feedback control (see \cite{GP11}).\\
\noindent
The numerical approximation of HJB equations has been studied with different techniques such as Finite Difference, Finite Volume and Semi-Lagrangian schemes. We refer the interested reader to \cite{FFbook} for a comprehensive analysis of these methods.\\
The DPP is known to suffer the so called {\em curse of dimensionality}, namely the computational complexity of the problem increases exponentially when the dimension does. The problem is much harder when dealing with PDEs since their spatial discretization leads to huge systems of ODEs. Typically, we are able to solve a HJB equations numerically up to dimension 4 or 5. For this reason, model reduction plays a crucial role in order to reduce the complexity of the problem and to make the control problems feasible.
The procedure is thus split in two parts, where the first part consists of finding a reduced order model (ROM) which is suitable for the control purpose, followed by the numerical solution of the HJB equations, associated with the control problem, where the full system is replaced with the ROM.\\
Proper Orthogonal Decomposition (POD, see \cite{Vol11}) and Balanced Truncation (BT, see \cite{A05}) are two of the most popular techniques for model reduction of dynamical systems, including spatially discretized PDEs. POD is a rather general method, which is based on a Galerkin projection method for nonlinear dynamical systems where the basis functions are built upon information on the system whereas the BT method is based on a Petrov-Galerkin projection, where the basis functions are obtained by solving two Lyapunov equations.
The latter approach is only valid for linear systems, although extensions can be formulated (see \cite{S93}).\\
The coupling between HJB equations and POD has already been proposed by a series of pioneering work \cite{KX05, KVX04}. A study of the feedback control and an adaptive method can be found in \cite{AF12} and \cite{AFK16}. Error estimation for the method has been recently studied in \cite{AFV15}. We refer to \cite{KK14} for the coupling with BT.\\
In addition to POD and BT, in this work we consider a new approach based on solutions of algebraic Riccati equations (ARE) for the approximation of the value function for linear quadratic problems. This approach turns out to better capture information of the control problem and improve the quality of the suboptimal control.
We analyze and compare the reduction techniques for linear and nonlinear dynamical systems. We note that in the nonlinear settings we linearize the dynamical system in a neighborhood of the desired state to apply BT and the MOR approach based on the solutions of the ARE equation.\\
The paper is organized as follows. In Section 2 we recall the main results on dynamic programming. Section 3 explains the model order reduction approaches and their application to the dynamic programming principle and the HJB equations. Finally, numerical tests are presented in Section 4 and conclusions are drawn in Section 5.

\section{Numerical approximation of HJB equations}
\label{sec:controlproblem}
In this section we recall the basic results for the approximation of the Bellman equation, more details can be found in  \cite{BC97} and \cite{FFbook}.\\
Let the dynamics be given by
\begin{eqnarray}\label{eq:dynconaut}
	\left\{\begin{array}{ll}
	\dot{y}(t)&=  f(y(t),u(t)),\quad t\geq 0,\cr
	y(0)&=x,
	 \end{array}\right.
\end{eqnarray}
where the state $y(t)\in \R^n$, the control $u(t) \in \R^m$ and $u\in\U\equiv\{ u:[0,+\infty)\rightarrow U,\,\text{measurable}\}$ where $U$ is a closed bounded subset of $\R^m$, and $x\in\R^n$ is the initial condition. If $f$ is Lipschitz continuous with respect to the state variable and continuous with respect to $(y,u)$, the classical assumptions for the existence and uniqueness result for the Cauchy problem \eqref{eq:dynconaut} are satisfied (see \cite{BC97}).\\
The cost functional $J:\U\rightarrow\R$ we want to minimize is given by:
\begin{eqnarray}\label{eq:joi}
	J_x(u(\cdot))\coloneqq\int_0^{\infty} g(y(s),u(s))\E^{-\lambda s} \D s\,,
\end{eqnarray}
where $g$ is Lipschitz continuous in both arguments and $\lambda \geq 0$ is a given parameter. The function $g$ represents the running costs and $\lambda$ is the discount factor which guarantees that the integral is finite whenever $g$ is bounded and $\lambda >0$. 
 Let us define the value function of the problem as
\begin{equation}
v(x)\coloneqq\inf_{u(\cdot)\in\U }J_x(u(\cdot))\,.
\end{equation}
The Dynamic Programming Principle (DPP) characterizes the value function as follows
\begin{equation}\label{dpp}
v(x)=\inf_{u\in\U}\{\int_0^T g(y_x(t,u),u(t))\E^{-\lambda t} \, \D t+v(y_x(T,u))\E^{-\lambda T}\},
\end{equation}
where $y_x(t,u)$ is the solution of the dynamics for a given initial condition $x$ and any $T>0$. From the DPP, one can obtain a characterization of the value function in terms of the following first order nonlinear Bellman equation
\begin{equation}\label{eq:hjb}
\lambda v(x)+\max_{u\in U} \{-f(x,u)\cdot Dv(x)- g(x,u)\}=0, \quad \hbox{ for }x\in\R^n\,.
\end{equation}
Here, $Dv(x)$ denotes the gradient of $v$ at the point $x$. Once the value function is computed we are able to build the feedback as follows:
$$u^*(x):=\argmin_{u\in U} \{f(x,u)\cdot Dv(x) + g(x,u)\}.$$
Several approximation schemes on a fixed grid $G$ have been proposed for \eqref{eq:hjb}. Here we will use a semi-Lagrangian approximation based on the Dynamic Programming Principle. This leads to
\begin{equation}\label{hjbh}
v_{\Delta t}(x) =\min_{u\in U} \{ \E^{-\lambda \Delta t} v_{\Delta t}\left(x+\Delta t f\left(x,u\right)\right)+ \Delta t g\left(x,u\right)\}\,,
\end{equation}
where $v_{\Delta t}(x)$ converges to $v(x)$ when $\Delta t\rightarrow 0$.
A natural way to solve \eqref{hjbh} is to write it in fixed point iteration form
\begin{equation}\label{hjbhk}
V_i^{k+1}=\min_{u\in U} \{ \E^{-\lambda \Delta t} \mathcal{I}[V^k]\left(x_i+\Delta t f\left(x_i,u\right)\right)+ \Delta t g\left(x_i,u\right)\}\,,\quad i=1,\ldots,N_G.
\end{equation}
\noindent Here $V^k_i$ represents the values of the value function $v$ at a node $x_i$ of the grid at the $k$-th iteration in \eqref{hjbhk} and $\mathcal{I}$ is a multilinear interpolation operator acting on the values of the equidistant grid $G$ with mesh spacing denoted by $\Delta x$.

The method is referred to in the literature as the {\em value iteration method}. The convergence of the value iteration can be very slow and accelerated techniques, such as the {\em policy iteration} technique, can be found in \cite{AFK15}.

\begin{remark}
  \label{rmrk:lqr}
  Let us mention that in general it is hard to find an explicit
  solution for equation \eqref{eq:hjb} due to the nonlinearity of the
  problem. A particular case is the so called linear quadratic
  regulator (LQR) problem where the dynamics is linear and the cost
  functional is quadratic. The equations are thus given as
$$f(y,u)=Ay+Bu, \quad g(y,u)=y^TQy+ u^TRu,$$
where $A,Q\in\R^{n\times n}$,
$B\in\R^{n\times m}$, $R\in\R^{m\times m}$ with $Q$ and $R$ symmetric and $Q$ positive semi-definite and $R$ positive definite. Furthermore,
the set of admissible control values is $U = \R$.
Under these assumptions, it is known that the value function at any point $x\in \R^n$ is given by $v(x)=x^T P x$
where $P\in \R^{n\times n}$ is the solution of the following shifted algebraic Riccati equation (ARE):
\begin{equation}\label{eq:ricc}
(A-\lambda I_n)^TP+P(A-\lambda I_n)-PBR^{-1}B^TP+Q=0.
\end{equation}
Here, $I_n\in\R^{n\times n}$ is the $n$-dimensional identity matrix.
Finally, the optimal control is directly given in an appropriate
state-feedback form $u(t) = -R^{-1}B^TPy(t)$.
More details on the LQR can be found in \cite{CZ95}. We will use the
LQR problem for comparison purposes as a benchmark model in the
numerical examples, see Section~4.
\end{remark}

\section{Model Reduction}
\label{sec:mor}
The focus of this section is to recall some model reduction
techniques utilized to build surrogate models in this work. The Reduced Order Modelling
(ROM) approach to optimal control problems is based on projecting the
nonlinear dynamics onto a low dimensional manifold utilizing
projectors that contain information of the expected controlled
dynamics. The idea behind the projection techniques is to restrict the
dynamics $y(t)$
onto a low-dimensional subspace $\cV \subset \R^n$
that contains the relevant information. We equip the space $\cV$
with a basis matrix $V \in \R^{n\times \ell}$,
and approximate the full state vector by $y(t) \approx V \yl(t)$,
where $\yl(t):[0,\infty)\rightarrow \R^\ell$
are the reduced coordinates. Plugging this ansatz into the dynamical
system~\eqref{eq:dynconaut}, and requiring a so called Petrov-Galerkin
condition yields
\begin{empheq}[left=\empheqlbrace]{gather}
  \label{eq:dynconautred}
  \begin{aligned}
  \dot{y}^\ell(t)&= W^Tf(V \yl(t),u(t))\cr
  \yl(0)&=W^Tx,
  \end{aligned}
\end{empheq}
where the matrix $W \in \R^{n \times \ell}$
is chosen, such that $W^TV=I_\ell$.
Further sampling based techniques can be employed to obtain an
efficient scheme for nonlinear problems as suggested in \cite{CS10, DHO12} and the reference therein.  
The presented procedure is a generic framework for model reduction.  It is
clear, that the quality of the approximation greatly depends on the
reduced space $\cV$.
In the next subsections, we briefly revisit some classical projection
techniques and introduce a new approach, which is tailored for the
approximation of the value function.

\subsection{Proper Orthogonal Decomposition}
\label{sec:mor_pod}
A common approach is based on the snapshot form of POD proposed
in \cite{Sir87}, which in the present situation works as
follows. We compute a set of snapshots $y_1,\dots,y_k$
of the dynamical system \eqref{eq:dynconaut} corresponding to
a prescribed input and different time instances $t_1,\ldots,t_k$ 
and define the POD ansatz of order $\ell$ for the state $y(t)$ by
\begin{equation}\label{pod_ans}
  y(t)\approx\sum_{i=1}^\ell y^\ell_i(t)\psi_i,
\end{equation}
where the basis vectors $\{\psi_i\}_{i=1}^\ell$
are obtained from the singular value decomposition of the snapshot
matrix $Y=[y_1,\ldots,y_k],$
i.e. $Y=\Psi\Sigma \Gamma$,
and the first $\ell$
columns of $\Psi$
form the POD basis functions of rank $\ell$.
Here the SVD is based on the Euclidean inner product. This is
reasonable in our situation, since the numerical computations
performed in our examples are based on Finite Difference
schemes.\\
In the present work the quality of the resulting basis is strongly related
to the choice of a given input $u$, 
whose optimal choice is usually unknown. 
For control problems, one way to
improve this selection is to compute snapshots from the following
equation for a given pair $(y,u)$ and any final time $T>0$
\begin{equation}\label{adj}
  -\dot{p}(t)=f_y(y(t),u(t))p(t)+g_y(y(t),u(t)),\quad p(T)=0,
 \end{equation}
 as suggested in \cite{SV13}. We refer to $p:[0,T]\rightarrow \R^n$ as the adjoint solution (see \cite{HPUU09}). The advantage of this approach is that
 it is able to capture the dynamics of the adjoint equation~\eqref{adj} which is
 directly related to the optimality conditions.

\subsection{Balanced truncation}
\label{sec:mor_bt}
The balanced truncation (BT) method is a well-established ROM
technique for LTI systems
\begin{align*}
\dot y(t) &= Ay(t)+Bu(t), \\
 z(t)&=Cy(t),
\end{align*}
where $z(t)$
is the output of interest. We refer to \cite{A05} for a complete
description of the topic.  The BT method is based on the solution of
the reachability Gramian $\tilde{P}$
and the observability Gramian $\tilde{Q}$
which solve respectively the following Lyapunov equations
$$ A\tilde{P} + \tilde{P} A^T +  B B^T = 0,\quad  A^T\tilde{Q}+\tilde{Q} A +  C^T C = 0.$$
We determine the Cholesky factorization of the Gramians
$$\tilde{P}=\Phi\Phi^T\qquad \tilde{Q}=\Upsilon\Upsilon^T.$$
Then, we compute the singular value decomposition of the Hankel operator $\Upsilon^T\Phi$ and set
$$W=\Upsilon U_1\Sigma_1^{1/2},\qquad V=\Upsilon V_1\Sigma_1^{1/2},$$
where $U_1, V_1\in\R^{n\times \ell}$
are the first $\ell$
columns of the left and right singular vectors of the Hankel operator
and $\Sigma_1=\mbox{diag}(\sigma_1,\ldots,\sigma_\ell)$
matrix of the first $\ell$ singular values.

The idea of BT is to neglect states that are
both, hard to reach and hard to observe. This is done by neglecting
states that correspond to low Hankel singular values $\sigma_i$.
This method is very popular, also because the whole procedure can 
be verified by a-priori error bounds in several system norms, and
the Lyapunov equations can be solved very efficiently due to their
typical low-rank structure in large-scale applications, see~\cite{BS13}.

\subsection{A new approach based on algebraic Riccati equations}
\label{sec:mor_ricc}
For arbitrary control problems, the value function is in general not
available in analytical form. However, in the case of LQR problem,
the value function has the quadratic form $v(x) = x^T P x$
where $P$
solves an algebraic Riccati equation~\eqref{eq:ricc}. 

Thus, the relevant information of the value function is stored in the
matrix $P$
and can be extracted by taking the SVD (or eigenvalue decomposition,
since $P$
is symmetric) $P = \Psi \Sigma \Psi^T$
with an orthonormal matrix $\Psi = [\psi_1, \dots, \psi_n]$.
We can approximate $P$
with $P^\ell = \sum_{k=1}^\ell \sigma_k \psi_k\psi_k^T$
and the error bound reads
$$ \| P - P^\ell \|_2 \leq \sigma_{\ell+1}, $$
where we applied the Eckart-Young-Mirsky theorem as mentioned in \cite{A05,SH16}.
We define the reduced value function as
$v^\ell(x) \coloneqq x^T P^\ell x$. 
Then the following bounds hold true
$$|v(x) - v^\ell(x)| \leq \sigma_{\ell+1}\|x\|^2, \quad \forall x \in \R^n.$$
Thus, if we define the reduced space
$\mathcal V \coloneqq \operatorname{span}(\psi_1, \dots, \psi_\ell)$,
we can expect an accurate approximation of the relevant information in the
value function, at least in the case where the system dynamics are
linear.  Furthermore, we note that the value function is of the form
$v(x)=x^TPx$
only in the case, where the set of controls is $U = \R$.

\subsection{The coupling between HJB and model reduction}

Since the curse of dimensionality prohibits a direct solution of the
HJB equations in higher dimensions, we apply model reduction in the
first place, in order to obtain a small system for which the HJB
equation admits a computable solution. In the general projection
framework above, we define the following reduced HJB problem, which is
the optimal control problem for the projected system:
\begin{gather}\label{KPl}
  \inf_{u \in \mathcal U} J^\ell_{W^T x} (u) = \inf_{u \in \mathcal U} \int_0^\infty g(V \yl(t),u(t),t)\E^{-\lambda t}\; \D t,\\
  \text{s.t.}\quad
\begin{aligned}
  \dot{\yl}(t) &= W^T f(V \yl(t), u(t)),\quad t \geq 0 \\
  \yl(0) &= W^Tx
\end{aligned}
\end{gather}
As in the full-dimensional case, the value function
$v^\ell(W^Tx) = \inf_{u \in \mathcal U} J^\ell_{W^Tx}(u)$
fullfills a $\ell$-dimensional HJB equations, which can be solved
numerically. This gives an approximation to the true (in general unknown)
value function at the point $x \in \R^n$:
\begin{align}
  \label{eq:vxapprox}
  \hat{v}^\ell(x) \coloneqq v^\ell(W^T x).
\end{align}
Furthermore, the reduced value function $\hat v^\ell(x)$
can be used to define a reduced feedback control function similar to
the full dimensional case as
$$ \hat u^*(x) \coloneqq \min_{u\in U} \{f(x,u)\cdot D\hat v^\ell(x) + g(x,u)\}. $$
\begin{remark}\label{rmrk:red_space}
  For the numerical approximation of the value function, we must restrict our computational domain in the $\ell$-dimensional reduced space. Since the physical meaning of the full-coordinates is
  lost when going to the reduced coordinates, it is in general not
  clear how to choose the interval lengths of the grid. We therefore
  restrict ourselves to the approximation of the value function for
  vectors in the set
  $\Theta \coloneqq \{x \in \R^n \mbox{ s.t. }  \|x\|_\infty \leq a\}$,
  i.e. for all $x\in\Theta$
  and $i=1,\dots,n$
  it holds $|x_i|\leq a$,
  where $x_i$
  denotes the $i$-th
  component of $x$.
  We then define the reduced domain, which is to be discretized as
  $\Theta_{\ell} \coloneqq \bigtimes_{i=1}^\ell (\underbar x_{i}, \bar x_{i})
  \subset \R^\ell$, where the interval boundaries $\underbar x_{i}$
  and $\bar x_{i}$
  are calculated in such a way that for all full states $x\in \Theta$,
  the projected vectors are mapped to vectors in $\Theta_\ell$,
  i.e. $W^T x \in \Theta_{\ell}$
  for all $x\in\Theta$.
  Thus, we expect to have a valid value function for all vectors
  $x\in\Theta$. A different approach for the reduced interval can be found in \cite{AF12}.
\end{remark}

\section{Numerical Examples}
\label{sec:errors}
We now compare the different approaches introduced in
Section~\ref{sec:mor}. The first example is a classical LQR scenario,
i.e. a linear system with quadratic cost functional. This simple
setup has the huge advantage of a known value function, that can be
used for comparing the different approaches for the
HJB approximations. 
In the second example, we study the behavior of
the feedback control for a nonlinear viscous Burgers equation.

\subsection{One-Dimensional Heat Advection-Diffusion Equation}
\label{sec:num1}
\begin{figure}[tb] 
  \begin{minipage}{.49\textwidth}
    \leftfigure{
      \includegraphics{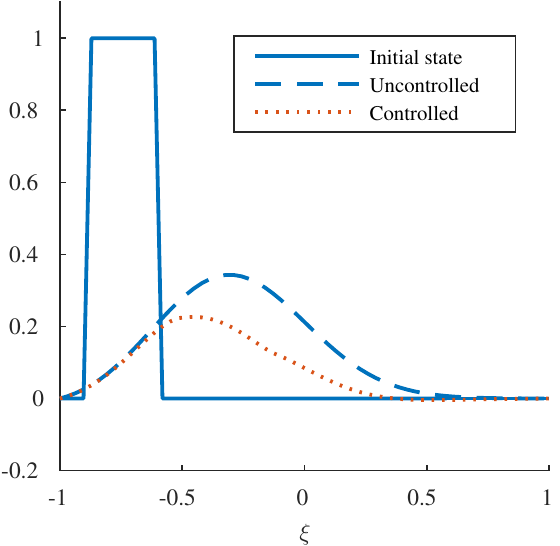}
    }
  \end{minipage}
  \hfill
  \begin{minipage}{.49\textwidth}
    \rightfigure{
     \includegraphics{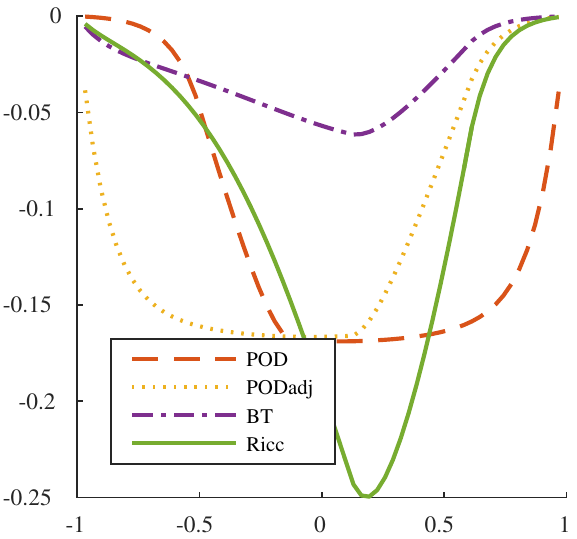}
    }
  \end{minipage}
  \leftcaption{Initial state, uncontrolled and controlled state of the
    linear example at time $t=0.3$.}
  \label{fig:lin_state}
  \rightcaption{Dominant basis vectors for all approaches. }
  \label{fig:linear_bases}
\end{figure}

Our first example consists of a one-dimensional advection-diffusion
equation 
\begin{align*}
  \partial_t w(t,\xi) - \mu_\text{diff} \partial_{\xi\xi} w(t,\xi) + \mu_\text{adv} \partial_\xi w(t,\xi) &= \mathbf 1_{\Omega_B}(\xi) u(t), &&t \geq 0, \xi\in \Omega \\
  w(t,\xi) &= 0, &&t \geq 0, \xi\in\{-1,1\} \\
  w(0,x) &= w_0(x),  &&x \in \Omega \\
  z(t) &= \frac{1}{|\Omega_C|}\int_{\Omega_C} w(t,\xi) \D \xi,
\end{align*}
\noindent
with $\Omega \coloneqq (-1, 1)\subset \R$
and distributed control acting on a set $\Omega_B =
[-0.5,-0.1]$. The output of interest $z(t)$
is the average temperature distribution on the interval
$\Omega_C = [0.1,0.6]$, $\mathbf 1_{\Omega_B}(\xi)$
and $\mathbf 1_{\Omega_C}(\xi)$
denote the characteristic functions of the set $\Omega_B$
resp. $\Omega_C$ at the point $\xi \in \Omega$.
We choose the parameter values $\mu_\text{diff}=0.2$
and $\mu_\text{adv}=2$.
We discretize the PDE in space by using a finite difference scheme on
an equidistant grid with interior points $\{\xi_i\}_{i=1}^n$. The dimension of the semi-discrete problem is $61$. The advection term is discretized by using an upwind scheme. In order to solve the problem numerically for the simulation and
the generation of the snapshots, we apply an explicit
Euler scheme. In order to obtain a control problem, we introduce the cost functional
as in Remark~\ref{rmrk:lqr} with $Q = 20C^TC$ and $R=0.1$, where $C$ is the
discretized representation of $z(\cdot)$. The final
setting is given by
\begin{gather*}
  \min_{u \in L^2(0,\infty)} \int_0^\infty (20 z(t)^2 + 0.1 u(t)^2) \D t\\
  \text{s.t.}\quad \dot y(t) = f(y(t),u(t)) = A y(t) + B u(t),\quad
  z(t)= C y(t), \quad y(0) = x.
\end{gather*}
The solution to this problem can be calculated in a closed loop form
and is given by $u(t) = -10 B^T P x(t)$,
where $P \in \R^{n\times n}$
solves the associated ARE~\eqref{eq:ricc} with $\lambda =
0$. Furthermore, the value function is known to be a quadratic
function of the form $v(x) = x^T P x$.
Figure~\ref{fig:lin_state} shows the controlled and uncontrolled
solution for the initial condition
$x = \left( 0.2 \cdot \mathbf 1_{(-0.8,-0.6)}(\xi_i)\right)_{i=1}^n$,
where the true LQR control is used to generate the figure.

We now construct the bases $W_q$
and $V_q$
for the different approaches $q\in\{$POD,
PODadj, BT, Ricc$\}$
introduced in Section~\ref{sec:mor}. In order to obtain the
basis for the POD approach, we simulate the full system with a
prescribed control function $u(t) = \sin(t)$
for $t\in[0,2\pi]$
and compute the POD method as explained in Section 3.  Since $\Wpod$
is an orthonormal matrix, we simply set the biorthogonal counterpart
as $\Vpod \coloneqq \Wpod$.
The basis $\Wpodadj$
for the adjoint system are calculated with the same control input and
discretization parameters, but solving equation \eqref{adj}. The basis
matrices for balanced truncation are denoted as $\Wbt$
and $\Vbt$
and are calculated in the usual way as explained in
Section~\ref{sec:mor}. Finally, the Riccati basis is built by taking
the first $\ell$
left singular vectors of the SVD of $P \in \R^{n\times n}$,
where $P$ solves the ARE~\eqref{eq:ricc}.

We now calculate the reduced value
functions for the different approaches, which we will denote as
$\hat v^\ell_{q}$
with $q$
as above. We apply a value iteration scheme based on an equidistant
grid in $\ell$
dimensions. For details, we refer to Section~\ref{sec:controlproblem}
and the references given there. The goal in this linear example is to
reproduce the true LQR control and value function by the HJB
approach.
The set of admissible controls is chosen as a discrete grid on the
interval $[-2, 2]$
with $301$
grid points. This set of controls is sufficiently large, to capture
the control values for all possible vectors $x \in \Theta$ with $a=0.2$, see Remark~\ref{rmrk:red_space}.

As a first qualitative comparison, we plot the dominant basis vectors
of all different approaches in Fig.~\ref{fig:linear_bases}. It can be
seen that the basis vectors carry very different
information. Especially the basis vector for the Riccati approach does
not reflect the input region of
the model very well, but it provides details about the region of
measurement $\Omega_C$. 
Still, by its construction we expect accurate approximations of the value function.

\begin{figure}
  \centering
  \begin{minipage}{\textwidth}
    \includegraphics{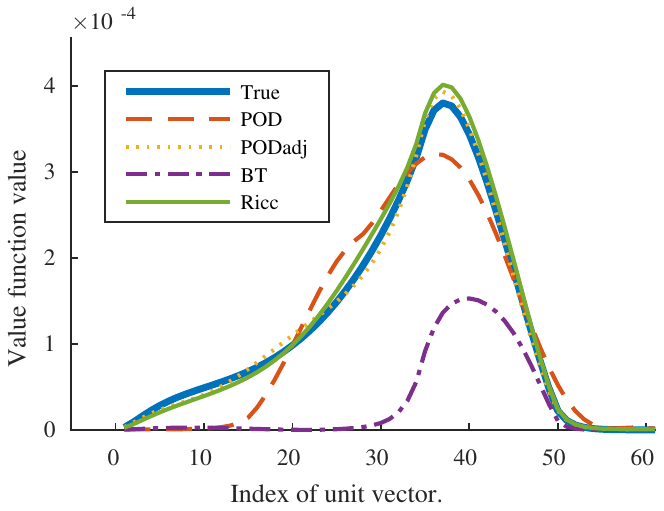}
    \includegraphics{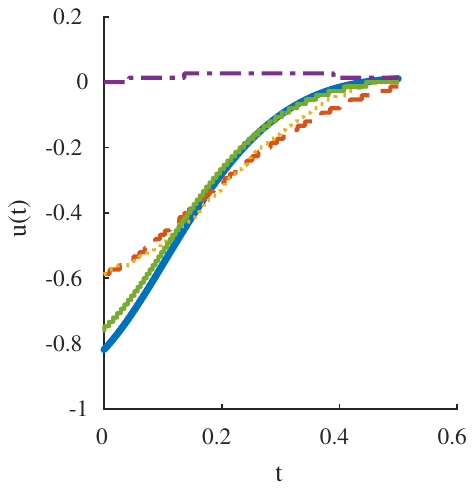}
  \end{minipage}
  \caption{Results for the approximation of the value function (left)
    for $\ell = 3$.
    True (LQR) control and approximated controls (right) for a a given
    initial state $x=0.2(B+C^T)$.}
  \label{fig:linear_vf_results}
\end{figure}

Another interesting insight is given, when we compare the values of
the approximated value functions $\hat v_{\ell,q}(\cdot)$
at the points $x_i \coloneqq 0.2 e_i$,
where $e_i$
is the $i$-th
unit vector in $\R^n$.
The results are depicted in Fig.~\ref{fig:linear_vf_results} for
$\ell=3$.
We see that the different bases deliver different results: the
Riccati and adjoint approach capture the original behavior of the
value function. We note that if we increase the dimension of the surrogate, the results improve for
all approaches. In Fig. \ref{fig:linear_vf_results} we also show the resulting optimal control, and again we can see how the Riccati and adjoint approach are able to recover the true control signal.

A more quantitative comparison is given in
Table~\ref{table:linear_rand_vf}: We calculate the values of the true
value function and the reduced value functions for all approaches for
$50$
random test vectors from the set $\Theta$.
We next calculate the relative error between the approximation and the
true LQR value function and list the mean approximation error in
Table~\ref{table:linear_rand_vf}.  In this example, the POD-basis does
not yield accurate approximations to the true value function. Balanced
trunction requires more basis functions to capture the relevant
information for the value function. Only the adjoint approach and the
basis $\Wricc$ yield very accurate results.
\begin{table}
  \centering
    \begin{tabular}{ccccc}                       
      \toprule                                     
      & $\ell=1$ & $\ell=2$ & $\ell=3$ & $\ell=4$ \\                          
      \midrule                          
POD & 0.6426 & 0.5634 & 0.3297 & 0.3752 \\   
PODadj & 0.8144 & 0.4008 & 0.1036 & 0.0959 \\
BT & 0.9971 & 0.8271 & 0.7387 & 0.5848 \\    
Ricc & 0.5472 & 0.1363 & 0.0711 & 0.0566 \\             
      \bottomrule                                  
    \end{tabular}
    \caption{Approximation of the value function for the different approaches.}
        \label{table:linear_rand_vf}
\end{table} 

\subsection{Viscous Burgers Equation}
\label{sec:num2}
\begin{figure}[tb]
  \begin{minipage}{\textwidth}
    \leftfigure{
      \includegraphics{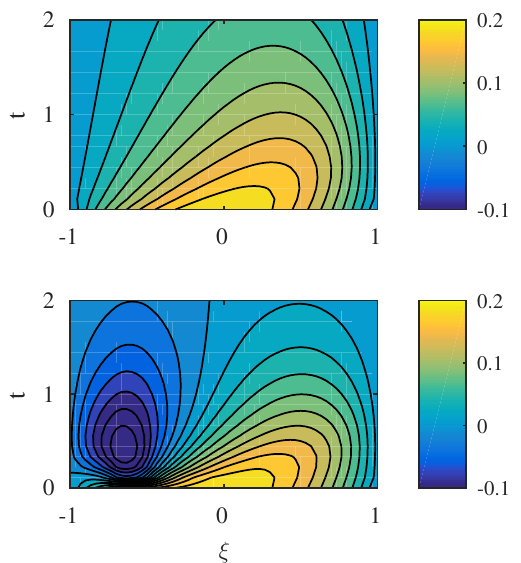}
    }
    \rightfigure{
      \includegraphics{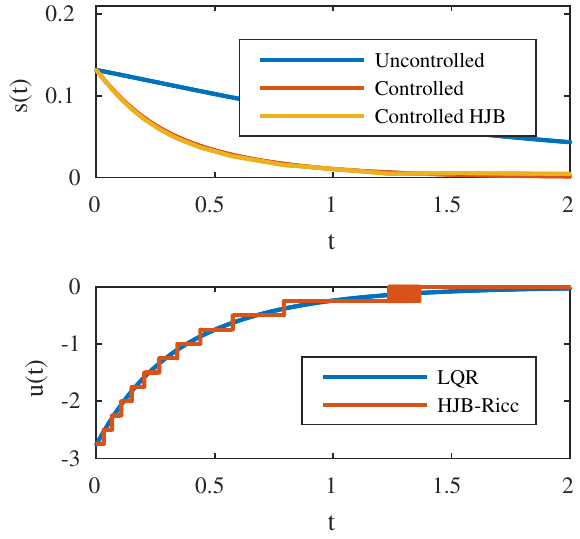}
    }
  \end{minipage}
  \leftcaption{Uncontrolled (top) and LQR-controlled (bottom) state
    example of the Burgers equation.}
  \label{fig:burgers_state}
  \rightcaption{Output of
    interest (top) and control (bottom) for the nonlinear Burgers
    example with the LQR and HJB-Riccati control for $\ell=4$.}
  \label{fig:burgers_output_control}
\end{figure}

Let us now study a more complex dynamical system, where no analytical value function can be derived. We choose
the 1D viscous Burgers equation on the domain
$\Omega \coloneqq (-1,1)$
with homogeneous Dirichlet boundary conditions. The continuous
equations now read as follows:
\begin{align*}
  \partial_t w(\xi,t) - 0.2 \partial_{\xi\xi} w(\xi, t) + 5 w \partial_\xi w(\xi,t) &= \mathbf 1_{\Omega_B}(\xi) u(\xi), \quad \xi \in \Omega, t\geq 0\\
  w(\xi,t) &= 0, \quad \xi \in\{-1,1\}, t\geq0,\\
  w(\xi,0) &= w_0(\xi),\quad \xi \in \Omega.
\end{align*}
The output of interest in this case is defined as the integral of the
state over the whole domain:
$z(t) \coloneqq \int_\Omega w(\xi,t) \D \xi$
for $t\geq 0$. The control acts on the subdomain $\Omega_B \coloneqq [-0.7, -0.5]$.
The semi-discretization is again performed by using finite
differences with the same setting as in the linear example.  The
discretized system has now dimension $n=61$
and all computations are again performed by using an explicit Euler
scheme. The discretized PDE and the discretized output then have the form~\eqref{eq:dynconaut} with
$$f(y(t),u(t)) = A y(t) + B u(t) + \tilde f(y(t)),\quad y(0) = x, \quad z(t) = C y(t),$$
where $\tilde f(y)$
models the discretized nonlinear transport term.\\
 We introduce an
infinite-horizon optimal control problem, similar to the LQR case, by
defining the cost functional for the discretized equations as 
$$ J_x(u(\cdot)) \coloneqq \int_0^\infty ( 100 z(t)^2 + 0.1 u(t)^2 )\E^{-\lambda t} \D t$$
with the discount factor $\lambda = 1$.
Figure~\ref{fig:burgers_state} shows the uncontrolled state and
controlled solution. We note that the stabilization of the Burgers
equation via LQR problem has been studied in \cite{BK91}. The control
in the latter case has been computed after a linearization of the
dynamics around the set point $y = 0$
in order to solve the ARE equation. The continuous initial condition
is $w_0(\xi) = 0.2(1 - \xi^2)$.
The corresponding output and control function is depicted in
Figure~\ref{fig:burgers_output_control}. We can observe that the
Riccati based approach is able to recover the LQR control. We recall
that in the HJB setting the control space is discretized and it is not
continuous as in the LQR setting.

We build the different bases for this example with the same setting as
in the linear example before, only the time-steps for the HJB scheme
have been adjusted and the controls are chosen as $41$ equidistant
points from $[-5,5]$ in order to allow the necessary higher
control values. For the BT and the Riccati approach, we linearize
the system around $y=0$
and obtain a heat equation for which the BT basis and the ARE solution
are calculated. Then, the calculation of the value
function is performed for the nonlinear reduced equation.

In this example we do not have a closed-loop form of the value
function and thus we need a different way to compare the results. For
this purpose, we approximate the value of the cost functional
numerically by performing a highly-resolved simulation, followed by a
quadrature using the trapezoidal rule. We simulate the closed-loop
systems until $T=5$, which suffices to neglect the increment in the
cost functional on $t\in(5,\infty)$.

To compare the methods we show in Table
\ref{tbl:burgers_prescribed_x0} the evaluation of the cost functional
for different initial conditions
$x_{0,1} = 0.2 B$ and $x_{0,2} = 0.2(1-\xi)^2$
 and model reduction methods. It is hard
to compare the method since we do not know the full solution, however
it turns out that the Riccati and POD adjoint approach have the minimum values
and are the closest to the full dimensional Riccati linearized control.

\begin{table}
  \centering
  \begin{tabular}{c|cccc|cccc}
\toprule
   & \multicolumn{4}{c}{$x_{0,1}$} &  \multicolumn{4}{c}{$x_{0,2}$}  \\
     & $\ell=1$ & $\ell=2$ & $\ell=3$ & $\ell=4$ & $\ell=1$ & $\ell=2$ & $\ell=3$ & $\ell=4$ \\
\midrule                                                                                   
Ricc & 0.2962 & 0.2958 & 0.2955 & 0.2956 & 0.3789 & 0.3786 & 0.3786 & 0.3785 \\            
POD & 0.3926 & 0.3171 & 0.3112 & 0.3006 & 0.4197 & 0.3817 & 0.3802 & 0.3790 \\             
BT & 0.2981 & 0.3169 & 0.3297 & 0.3260 & 0.3785 & 0.3987 & 0.4115 & 0.4080 \\              
PODadj & 0.2960 & 0.2958 & 0.2955 & 0.2953 & 0.3786 & 0.3786 & 0.3786 & 0.3786 \\          
LQR & 0.2959 & 0.2959 & 0.2959 & 0.2959 & 0.3786 & 0.3786 & 0.3786 & 0.3786 \\                \bottomrule                       
  \end{tabular}
  \caption{Cost functional values for different inital vectors.}
  \label{tbl:burgers_prescribed_x0}
\end{table}

\section{Conclusion}
\label{sec:conc}
In this paper we propose a comparison of different model order
reduction techniques for dynamic programming equations.  Numerical
experiments show that the POD adjoint and the Riccati based approach
provide very accurate approximation for the control problem with
quadratic cost functional. This is what one can expect since both
methods contain information about the optimization problem, unlike BT
and POD when the snapshots are generated with a random initial input.
Moreover, the Riccati based approach can be generalized to nonlinear
dynamics. Here we propose to linearize the system around one point of
interests. In the future we would like to investigate a greedy
strategy to select more points. A parametric scenario will also be
considered in a future work as proposed in \cite{SH16} for linear
dynamical systems.

\begin{acknowledgement}
The first author is supported by US Department of Energy grant number
  DE-SC0009324. The second and third authors thank the Baden W\"urttemberg Stiftung
  gGmbH and the German Research Foundation (DFG) for financial support within
  the Cluster of Excellence in Simulation Technology (EXC 310/1) at the University
  of Stuttgart.
\end{acknowledgement}



\begin{thebibliography}{99}
%
\bibitem{AF12}
A. Alla and M. Falcone. {\em An adaptive POD approximation method for the control of advection-diffusion equations}, in K. Kunisch, K. Bredies, C. Clason, G. von Winckel, (eds) {\em Control and Optimization with PDE Constraints,} International Series of Numerical Mathematics, \textbf{164}, Birkh\"auser, Basel, 2013, 1-17.
%
%
%
\bibitem{AFK15}
A. Alla, M. Falcone, and D. Kalise. {\em An efficient policy iteration algorithm for dynamic programming equations},  SIAM J. Sci. Comput., \textbf{37}, 2015, 181-200. 

\bibitem{AFK16}
A. Alla, M. Falcone and D. Kalise. {\em A HJB-POD feedback synthesis approach for wave equation}, 
Bulletin of the Brazilian Mathematical Society, New Series, {\bf47}, 2016, 51-64.


\bibitem{AFV15}
A. Alla, M. Falcone and S. Volkwein. 
{\em Error Analysis for POD approximations of infinite horizon problems via the dynamic programming principle.} Submitted to SIAM Journal on Control and Optimization, 2015.
%
%
%
%
\bibitem{A05}
A.C. Antoulas. {\em Approximation of Large-Scale Dynamical Systems}, SIAM, 2005.
%
\bibitem{AK01}
J.A. Atwell and B.B. King. {\em Proper orthogonal decomposition for reduced basis feedback controllers for parabolic equations}, Mathl Comput. Modelling, \textbf{33}, 2001, 1-19. 
%
\bibitem{BC97}
M. Bardi and I. Capuzzo-Dolcetta. {\em Optimal Control and Viscosity Solutions of Hamilton-Jacobi-Bellman Equations}. Birkh\"auser, Basel, 1997.
%
%
%
\bibitem{BS13}
P. Benner and J. Saak. {\em Numerical solution of large and sparse continuous time algebraic matrix {R}iccati and {L}yapunov equations: a state of the art survey}, GAMM-Mitteilungen, 2013, 32-52.

\bibitem{BK91}
J. Burns and S. Kang. {\em A control problem for Burgers' equation with bounded input/output}, Nonlinear Dynamics {\bf 2}, 1991,235-262. 

\bibitem{CS10}
S. Chaturantabut and D.C. Sorensen. {\em Nonlinear model reduction via discrete empirical interpolation}, SIAM Journal on Scientific Computing, \textbf{32}, 2010, 2737-2764.
%
\bibitem{CZ95}
R.F. Curtain and H.J. Zwart. {\em An Introduction to Infinite-Dimensional Linear Systems Theory}, Springer, 1995.

\bibitem{DHO12}
M. Drohmann, B. Haasdonk and M. Ohlberger.
{\em Reduced Basis Approximation for Nonlinear Parametrized Evolution Equations based on Empirical Operator Interpolation},  SIAM J. Sci. Comput., {\bf 34}, 2012, 937-969.


%
%
%
\bibitem{FFbook}
M. Falcone and R. Ferretti. {\em Semi-Lagrangian Approximation Schemes for Linear and Hamilton-Jacobi equations}, 
SIAM, 2014.  

\bibitem{GP11}
L. Gr\"une, J. Panneck. {\em Nonlinear Model Predictive Control: Theory and Applications}, Springer, 2011.

%
\bibitem{HLBR12}
P. Holmes, J.L. Lumley, G. Berkooz, and C.W. Rowley. {\em Turbulence, Coherent Structures, Dynamical Systems and Symmetry}, Cambridge Monographs on Mechanics, Cambridge University Press, second edition, 2012.
%
\bibitem{HPUU09} 
M. Hinze, R. Pinnau, M. Ulbrich, and S. Ulbrich. {\em Optimization with PDE Constraints. Mathematical Modelling: Theory and Applications}, \textbf{23}, Springer Verlag, 2009.
%
%
\bibitem{KK14}
D. Kalise and A. Kr\"oner. {\em Reduced-order minimum time control of advection-reaction-diffusion systems via dynamic programming}, In Proceedings of the 21st International Symposium on Mathematical Theory of Networks and Systems, 2014, 1196-1202. 
%
%
\bibitem{KVX04}
K. Kunisch, S. Volkwein, and L. Xie. {\em HJB-POD based feedback design for the optimal control of evolution problems}, SIAM J. on Applied Dynamical Systems, \textbf{4}, 2004, 701-722.
%
\bibitem{KX05}
K. Kunisch and L. Xie. {\em POD-based feedback control of Burgers equation by solving the evolutionary HJB equation,} Computers and Mathematics with Applications, \textbf{49}, 2005, 1113-1126.

\bibitem{S93}
J. Scherpen.
{\em Balancing for nonlinear systems,} Systems Control Lett., {\bf 21}, 1993, 143-153.

\bibitem{SH16}
A. Schmidt, B. Haasdonk.
{\em Reduced Basis Approximation of Large Scale Algebraic Riccati Equations,} Simetech Preprint, Univ Stuttgart, 2015.

\bibitem{Sir87} L. Sirovich, {\em Turbulence and the dynamics of coherent structures. Parts I-II,}
Quarterly of Applied Mathematics, {\bf XVL} (1987), 561-590.

\bibitem{SV13}
A. Studinger and S. Volkwein. {\em Numerical Analysis of POD A-Posteriori Error Estimation for Optimal Control}, in K. Kunisch, K. Bredies, C. Clason, G. von Winckel, (eds) {\em Control and Optimization with PDE Constraints,} International Series of Numerical Mathematics, \textbf{164}, Birkh\"auser, Basel, 2013, 137-158.


%
\bibitem{Tro10}
F.~Tr\"oltzsch. {\em Optimal Control of Partial Differential Equations: Theory, Methods and Application}, American Mathematical Society, 2010.

\bibitem{Vol11}
S. Volkwein. {\em Model Reduction using Proper Orthogonal Decomposition}, Lecture Notes, University of Konstanz, 2013.\hfill\\
\texttt{http://www.math.uni-konstanz.de/numerik/personen/volkwein/\\teaching/scripts.php}

\end{thebibliography}
\end{document}